\documentstyle[12pt,epsf]{article}
\newcommand{\rr}{{\rm{I \! R}}}

\newcommand{\al}{\alpha}
\newcommand{\ds}{\displaystyle}
\newcommand{\f}{\frac}
\newcommand{\lam}{\lambda}

\newcommand{\om}{\omega}

\newcommand{\ca}{{\cal A}}
\newcommand{\th}{\theta}

\title{\bf Hopf bifurcation analysis for a mathematical model of P53-MDM2 interaction}
\author{Mihaela Neam\c tu$^{a}$\thanks{Corresponding author}, Raul Florin Horhat$^{b}$, Dumitru Opri\c s$^c$}
\date{}

\begin{document}
\maketitle

\begin{tabular}{cccccccc}
\scriptsize{$^{a}$Department of Economic Informatics, Mathematics and Statistics, Faculty of Economics,}\\
\scriptsize{West University of Timi\c soara, str. Pestalozzi, nr. 16A, 300115, Timi\c soara, Romania,}\\
\scriptsize{E-mail:mihaela.neamtu@fse.uvt.ro,}\\
\scriptsize{$^{b}$ Department of Biophysics and Medical Informatics,}\\
\scriptsize{University of Medicine and Pharmacy, Piata Eftimie Murgu, nr. 3, 300041, Timi\c soara, Romania,}\\
\scriptsize{E-mail: rhorhat@yahoo.com}\\
\scriptsize{$^{c}$ Department of Applied
Mathematics, Faculty of Mathematics,}\\
\scriptsize{West University of Timi\c soara, Bd. V. Parvan, nr. 4, 300223, Timi\c soara, Romania,}\\
\scriptsize{E-mail: opris@math.uvt.ro}\\

\end{tabular}

\bigskip

\begin{center}
\small{\bf Abstract}
\end{center}

\begin{quote} \small{In this paper we analyze a simple mathematical model which describes
the interaction between the proteins p53 and Mdm2. For the
sta\-tio\-na\-ry state we discuss the local stability and the
existence of the Hopf bifurcation. Choosing the delay as a
bifurcation parameter we study the direction and stability of the
bifurcating periodic solutions. Some numerical examples and the
conclusions are finally made.}
\end{quote}

\noindent{\small{{\it Keywords:} delay differential equation,
stability, Hopf bifurcation, p53, Mdm2.

\noindent{\it 2000 AMS Mathematics Subject Classification:34C23,
34C25, 37G05, 37G15, 92D10} .}}

\section*{\normalsize\bf 1. Introduction}

\hspace{0.6cm} In every normal cell, there is a protective
mechanism against tumoral degeneration. This mechanism is based on
the P53 network. P53, also known as "the guardian of the genome",
is a gene that codes a protein in order to regulate the cell
cycle. The name is due to its molecular mass: it is a 53 kilo
Dalton fraction of cell proteins. Mdm2 gene plays a very important
role in P53 network. It regulates the levels of intracellular P53
protein concentration through a feedback loop. Under normal
conditions the P53 levels are kept very low. When there is DNA
damage the levels of P53 protein rise and if there is a prolonged
elevation the cell shifts to apoptosis, and if there is only a
short elevation the cell cycle is arrested and the repair process
is begun. The first pathway protects the cell from tumoral
transformation when there is a massive DNA damage that cannot be
repaired, and the second pathway protects a number of important
cells (neurons, myocardic cells) from death after DNA damage. In
these cells first pathway, apoptosis, is not an option because
they do not divide in adult life and their importance is obvious.
Due to its major implication in cancer prevention and due to the
actions described above, P53 has been intensively studied in the
last two decades.

During the years, several models which describe the interaction
between P53 and Mdm2 have been studied. We mention some of them in
references [2], [3], [4], [6], [7], [8], [9], [10]. This paper
gives a mathematical approach to the model described in [10]. The
authors of paper [10] make a molecular energy calculation based on
the classical force fields, and they also use chemical reactions
constants from literature. Their results obtained by simulations
in accordance with experimental behavior of the P53-Mdm2 complex,
but they lack the mathematical part which we develop in this
paper. We analyze the Hopf bifurcation with time delay as a
bifurcation parameter using the methods from [1], [5], [7].

The paper is organized as follows. We present the mathematical
model in section 2 and the existence of the stationary state is
study. In section 3, we discuss the local stability for the
stationary state of system (2) and we investigate the existence of
the Hopf bifurcation for system (2) using time delay as the
bifurcation parameter. In section 4, the direction of Hopf
bifurcation is analyzed according to the normal form theory and
the center manifold theorem introduced by Hassard [5]. Numerical
simulations for confirming the theoretical results are illustrated
in section 5. Finally, some conclusions are made.

\section*{\normalsize\bf 2. The mathematical model and the stationary state}

\hspace{0.6cm}

The state variables are: $y_1(t), y_2(t)$ the total number of p53
molecules and the total number of Mdm2 proteins.

The interaction function between P53 and MDM2 is
$f:\rr^2_+\rightarrow\rr$ given by [10]:
$$f(y_1, y_2)=\ds\f{1}{2}(y_1+y_2+k-\sqrt
{(y_1+y_2+k)^2-4y_1y_2}).\eqno(1)$$

The parameters of the model are: $s$ the production rate of p53,
$a$ the degradation rate of p53 (through ubiquitin pathway), and
also the rate at which Mdm2 re-enters the loop, $b$ the
spontaneous decay rate of p53, $d$ the decay rate of the protein
rate Mdm2, $k_1$ the dissociation constant of the complex
P53-Mdm2, $c$ the constant of proportionality of the production
rate of Mdm2 gene with the probability that the complex P53-Mdm2
is build. These parameters are positive numbers.

The mathematical model is described by the following differential
system with time delay [10]:

$$\begin{array}{l}
\vspace{0.1cm}
\dot y{}_1(t)=s-af(y_1(t), y_2(t))-by_1(t),\\
\dot y{}_2(t)=cg(y_1(t-\tau),
y_2(t-\tau))-dy_2(t),\end{array}\eqno(2)$$ where $f$ is given by
(1) and $g:\rr^2_+\rightarrow\rr$, is
$$g(y_1, y_2)=\ds\f{y_1-f(y_1, y_2)}{k_1+y_1-f(y_1,y_2)}.\eqno(3)$$

For the study of the model (1) we consider the following initial
values:

$$y_1(\th)=\varphi_1(\th), y_2(\th)=\varphi_2(\th),\th\in[-\tau,0],
$$ with $\varphi_1, \varphi_2:[-\tau, 0]\rightarrow\rr_+$ are differentiable functions.

In the second equation of (2) there is delay, because the
trancription and translation of Mdm2 last for some time after that
p53 was bound to the gene.

The stationary state $(y_{10}, y_{20})\in \rr^2_+$ is given by the
solution of the system of equations:

$$\begin{array}{l}
s-af(y_1, y_2)-by_1=0,\\
cg(y_1,y_2)-dy_2=0.\end{array}\eqno(4)$$

From (1), (3) and (4) we deduce that the stationary state can be
found through the intersection of the curves:
$$y_2=f_1(y_1), \quad y_1=f_2(y_1),\eqno(5)$$
where $f_1, f_2:\rr_+\rightarrow\rr$ are given by:
$$\begin{array}{l}
f_1(y_1)=\ds\f{c(y_1(a+b)-s)}{d(k_1a+(a+b)y_1-s)},\\
f_2(y_2)=\ds\f{(s-by_1)(ka+(a+b)y_1-s)}{a((a+b)y_1-s)}.\end{array}\eqno(6)$$

{\bf Proposition 1}. {\it The function $f_1, f_2$ from (6) have
the following properties:

(i) $f_1$ is strictly increasing, $f_2$ is strictly decreasing on
$\left (\ds\f{s}{a+b}, \ds\f{s}{b}\right )$;

(ii) There is a unique value $y_{10}\in \left (\ds\f{s}{a+b},
\ds\f{s}{b}\right )$ so that $f_1(y_{10})=f_2(y_{10})$, where
$y_{10}$ is the solution of the equation $\varphi(x)=0$, from
$\left (\ds\f{s}{a+b}, \ds\f{s}{b}\right )$ where}
$$\begin{array}{l}
\varphi(x)\!=\!ac((a\!+\!b)x\!-\!s)^2\!-\!d(k_1a\!+\!(a\!+\!b)x\!-\!s)(ka\!+\!(a\!+\!b)x\!-\!s)(s\!-\!bx).\end{array}\eqno(7)$$
{\bf Proof}. (i) From (6) we have:
$$\begin{array}{l}
f'_1(y_1)=\ds\f{(a+b)ack_1}{d(k_1a+(a+b)y_1-s)^2},\\
f'_2(y_2)=-\ds\f{b}{a}-\ds\f{bk}{(a+b)y_1-s}-\ds\f{k(a+b)(s-by_1)}{((a+b)y_1-s)^2}.\end{array}\eqno(8)$$
For $y_1\in \left (\ds\f{s}{a+b}, \ds\f{s}{b}\right )$, the
relations (8) lead to $f'_1(y_1)>0$, $f'_2(y_1)<0$, so that $f_1$
is strictly increasing and $f_2$ is strictly deceasing.

(ii) By (i) there is $y_{10}\in \left (\ds\f{s}{a+b},
\ds\f{s}{b}\right )$ so that $f_1(y_{10})=f_2(y_{10}).$ From (6)
and (7) it results that:
$$\varphi(y_1)=f_{12}(y_1)h(y_1),$$ where
$$\begin{array}{l}
h(y_1)=ad((a+b)y_1-s)((a+b)y_1+k_1a-s)\\
f_{12}(y_1)=f_1(y_1)-f_2(y_1).\end{array}$$ On $\left
(\ds\f{s}{a+b}, \ds\f{s}{b}\right )$, the functions $f_{12}$ and
$h$ are strictly increasing and consequently $\varphi$ is strictly
increasing. Because
$\varphi(\ds\f{s}{a+b})=-\ds\f{a^3dkk_1s}{a+b}<0$ and
$\varphi(\ds\f{s}{b})=\ds\f{a^3cs^2}{b^2}>0$ we can conclude that
equation $\varphi(x)=0$ has a unique solution $y_{10}$ on $\left
(\ds\f{s}{a+b}, \ds\f{s}{b}\right )$.

\section*{\normalsize\bf{3. The analysis of the stationary state and the existence of the Hopf bifurcation.}}

\hspace{0.6cm}

We consider the following translation:
$$y_1=x_1+y_{10}, y_2=x_2+y_{20}$$ and system (2) can be
expressed as:
$$\begin{array}{l}
\vspace{0.1cm}
\dot x{}_1(t)=s-af(x_1(t)+y_{10}, x_2(t)+y_{20})-b(y_1(t)+y_{10}),\\
\vspace{0.1cm}
\dot x{}_2(t)=cg(x_1(t-\tau)+y_{10}, x_2(t-\tau)+y_{20})-d(x_2(t)+y_{20}).\\
\end{array}\eqno(9)$$

System (9) has a unique stationary state $(0,0)$. To investigate
the local stability of the equilibrium state we linearize system
(9). Let $u_1(t)$ and $u_2(t)$ be the linearized system variables.
Then (9) is rewritten as:

$$\dot U(t)=AU(t)+BU(t-\tau),\eqno(10)$$
where
$$A\!\!=\!\!\left(\!\!\!\!\begin{array}{cc}
\vspace{0.2cm}
-(b+a\rho_{10}) & -a\rho_{01} \\
\vspace{0.2cm}
0 & -d\\
\end{array}\!\!\!\!\right)\!\!, B\!\!=\!\!\left(\!\!\!\!\begin{array}{cc}
\vspace{0.2cm}
0 & 0 \\
\vspace{0.2cm} c\gamma_{10} &
c\gamma_{01}\end{array}\!\!\!\!\right)\eqno(11)$$ with
$U(t)=(u_1(t), u_2(t))^T$ and $\rho_{10}$, $\rho_{01}$,
$\gamma_{10}$, $\gamma_{01}$ are the values of the first order
derivatives for the functions:
$$f(x,y)=\ds\f{1}{2}(x+y+k-\sqrt{(x+y+k)^2-4xy})$$
and $$g(x,y)=\ds\f{x-f(x,y)}{k_1+x-f(x,y)}$$ evaluated at
$(y_{10}, y_{20})$.

The characteristic equation corresponding to system (10) is
$\Delta(\lambda, \tau)=det(\lambda I-A-e^{-\lambda \tau}B)=0$
which leads to:
$$\lam^2+p_1\lam+p_0-(q_1\lam+q_0) e^{-\lam \tau}=0,\eqno(12)$$
where
$$\begin{array}{l}
p_1=b+d+a\rho_{10}, \quad p_0=d(b+a\rho_{10}), \quad q_1=c\gamma_{01},\\
q_0=c\gamma_{01}(b+a\rho_{10})-ac\rho_{01}\gamma_{10}
.\end{array}$$

If there is no delay, the characteristic equation (12) becomes:
$$\Delta(\lam, 0)=\lam^2+(p_1-q_1)\lam+p_0-q_0.\eqno(13)$$
Then, the stationary state $(0,0)$ is locally asymptotically
stable if
$$p_1-q_1>0, \quad p_0-q_0>0.\eqno(14)$$

When $\tau>0$, the stationary state is asymptotically stable if
and only if all roots of equation (12) have a negative real part.
We are determining the interval $[0, \tau_0)$ so that the
stationary state remain asymptotically stable.

In what follows we study the existence of the Hopf bifurcation for
equation (10) choosing $\tau$ as the bifurcation parameter. We are
looking for the values $\tau_0$ of $\tau$ so that the stationary
state $(0,0)$ changes from local asymptotic stability to
instability or vice versa. We need the pure imaginary solutions of
equation (12). Let $\lam=\pm i\om_0$ be these solutions and
without loss of generality we assume $\om_0>0$. Replacing
$\lam=i\om_0$ and $\tau=\tau_0$ in (12) we obtain:

$$\begin{array}{l}
q_0cos\om_0\tau_0+q_1\om_0sin\om_0\tau_0=p_0-\om_0^2\\
q_0sin\om_0\tau_0-q_1\om_0cos\om_0\tau_0=-\om_0p_1,\end{array}$$
which implies that
$$\tau_0\!=\!\ds\f{1}{\om_0}\!\left (\!(2k\!+\!1)\pi\!+\!arcsin\ds\f{p_1\om_0}
{\sqrt{(p_0\!-\!\om_0^2)^2\!+\!\om_0^2p_1^2}}\!+\!arcsin\ds\f{q_1\om_0}
{\sqrt{(p_0\!-\!\om_0^2)^2\!+\!\om_0^2p_1^2}}\!\right )\eqno(15)$$
where $\omega_0$ is a solution of the equation:
$$\om^4+(-p_1^2-2p_0+q_1^2)\om^2+p_0^2-q_0^2=0.$$

Now we have to calculate $Re\left (\ds\f{d\lam}{d\tau}\right )$
evaluated at $\lam=i\om_0$ and $\tau=\tau_0$. We have:
$$\ds\f{d\lam}{d\tau}|_{\lam=i\om_0,\tau=\tau_0}=M+iN$$ where
$$M=\ds\f{q_1^2\om_0^6+2q_0^2\om_0^4+(p_1^2q_0^2-p_0^2q_1^2-2p_0q_0^2)\om_0^2}
{l_1^2+l_2^2}\eqno(16)$$ and
$$\begin{array}{lll}
N & = &
\ds\f{-q_1^2\tau_0\om_0^7+\om_0^5(q_0q_1-p_1q_1^2\!+\!\tau_0(2p_0q_1^2-p_1^2q_1^2-q_0^2))}{l_1^2+l_2^2}+\\
& + &
\ds\f{\om_0^3(-p_1q_0^2\!-\!2p_0q_0q_1+p_1^2q_0q_1-p_0p_1q_1^2\!+\!\tau_0(\!-\!p_1^2q_0^2\!-\!q_1^2p_0^2\!+\!2p_0q_0^2))}{l_1^2+l_2^2}+\\
& + &
\ds\f{\om_0(-p_0p_1q_0^2+p_0^2q_0q_1\!-\!\tau_0p_0^2q_0^2)}{l_1^2+l_2^2}
\end{array}\eqno(17)$$
with
$$\begin{array}{l}
l_1=-q_1\om_0^2\!+\!p_1q_0\!-\!q_1p_0+\tau_0(-q_1p_1\om_0^2\!-\!q_0\om_0^2\!+\!q_0p_0),\\
l_2=2\om_0q_0\!+\!\tau_0(-q_1\om_0^3\!+\!p_0q_1\om_0\!+\!p_1q_0\om_0).\end{array}$$
We conclude with:

{\bf Theorem 1.} {\it If there is no delay, under condition (14)
system (10) has an asymtotically stable stationary state. If $\tau
>0$ and $p_1^2q_0^2-q_1^2p_0^2-2p_0q_0^2>0$ then there is
$\tau=\tau_0$ given by (15) so that $Re\left
(\ds\f{d\lam}{d\tau}\right )_{\lam=i\om_0, \tau=\tau_0}>0$ and
therefor a Hopf bifurcation occurs at $(y_{10}, y_{20})$.}

%\begin{center}
\section*{{\normalsize\bf 4. Direction and stability of the Hopf bifurcation}}
%\end{center}

In this section we describe the direction, stability and the
period of the bifurcating periodic solutions of system (2). The
method we use is based on the normal form theory and the center
manifold theorem introduced by Hassard  [5]. Taking into account
the previous section, if $\tau=\tau_0$ then all roots of equation
(11) other than $\pm i\om_0$  have negative real parts, and any
root of equation (11) of the form $\lam(\tau)=\al(\tau)\pm
i\om(\tau)$ satisfies $\al(\tau_0)=0$, $\om(\tau_0)=\om_0$ and
$\ds\f{d\al(\tau_0)}{d\tau}\neq0$. For notational convenience let
$\tau=\tau_0+\mu, \mu\in\rr$. Then $\mu=0$ is the Hopf bifurcation
value for equations (2).

The Taylor expansion at $(y_{10}, y_{20})$ of the right members
from (9) until the third order leads to:
$$\dot x(t)=Ax(t)+Bx(t-\tau)+F(x(t), x(t-\tau)),\eqno(18)$$ where
$x(t)=(x_1(t), x_2(t))^T$, $A$, $B$ are given by (11) and
$$F(x(t),x(t-\tau))=(F^1(x(t)),F^2(x(t-\tau)))^T,\eqno(19)$$ where
$$\begin{array}{lll}
F^1(x_1(t),x_2(t)) & = & -\ds\f{a}{2}[\rho_{20}x_1^2(t)+2\rho_{11}x_1(t)x_2(t)+\rho_{02}x_2^2(t)]-\\
& - &
\ds\f{a}{6}[\rho_{30}x_1^3(t)+3\rho_{21}x_1^2(t)x_2(t)+3\rho_{12}x_1(t)x_2^2(t)+\rho_{03}x_2^3(t)]\end{array}$$
$$\begin{array}{lll}
F^2(x_1(t-\tau),x_2(t-\tau)) & = & \ds\f{c}{2}[\gamma_{20}x_1^2(t-\tau)+2\gamma_{11}x_1(t-\tau)x_2(t-\tau)+\\
& + & \gamma_{02}x_2^2(t-\tau)]+\\
& + &
\ds\f{c}{6}[\gamma_{30}x_1^3(t-\tau)+3\gamma_{21}x_1^2(t-\tau)x_2(t-\tau)+\\
& + &
3\gamma_{12}x_1(t-\tau)x_2^2(t-\tau)+\gamma_{03}x_2^3(t-\tau)]\end{array}$$
and $\rho_{ij}$, $\gamma_{ij}$, $i,j=0,1,2,3$ are the values of
the second and third order derivatives for the functions $f(x,y)$
and $g(x,y)$.

 Define the space of
continuous real-valued functions as $C=C([-\tau_0,0],\rr^4).$

In $\tau=\tau_0+\mu, \mu\in\rr$, we regard $\mu$ as the
bifurcation parameter. For $\Phi\in C$ we define a linear
operator:
$$L(\mu)\Phi=A\Phi(0)+B\Phi(-\tau)$$ where A and B are given by
(11) and a nonlinear operator $F(\mu, \Phi)=F(\Phi(0)$,
$\Phi(-\tau))$, where $F(\Phi(0), \Phi(-\tau))$ is given by (19).
According to the Riesz representation theorem, there is a matrix
whose components are bounded variation functions, $\eta(\theta,
\mu)$ with $\theta\in[-\tau_0, 0]$ so that:
$$L(\mu)\Phi=\int\limits_{-\tau_0}^0d\eta(\theta,\mu)\phi(\theta),
\quad \theta\in[-\tau_0,0].$$

For $\Phi\in C^1([-\tau_0, 0], \rr^{4})$ we define:
$$\ca(\mu)\Phi(\th)=\left\{\begin{array}{ll} \vspace{0.2cm}
\ds\f{d\Phi(\th)}{d\th}, & \th\in[-\tau_0,0)\\
\int\limits_{-\tau_0}^0d\eta(t,\mu)\phi(t), &
\th=0,\end{array}\right.$$
$$R(\mu)\Phi(\th)=\left\{\begin{array}{ll} \vspace{0.2cm}
0, & \th\in[-\tau_0,0)\\
F(\mu, \Phi), & \th=0.\end{array}\right.$$

We can rewrite (18) in the following vector form:
$$\dot u_t=\ca(\mu)u_t+R(\mu)u_t\eqno(20)$$
where $u=(u_1, u_2)^T$, $u_t=u(t+\theta)$ for $\theta\in[-\tau_0,
0]$.

For $\Psi\in C^1([0,\tau_0], \rr^{*4})$, we define the adjunct
operator $\ca^*$ of $\ca$ by:
$$\ca^*\Psi(s)=\left\{\begin{array}{ll} \vspace{0.2cm}
-\ds\f{d\Psi(s)}{ds}, & s\in(0, \tau_0]\\
\int\limits_{-\tau_0}^0d\eta^T(t,0)\psi(-t), &
s=0.\end{array}\right.$$ We define the following bilinear form:
$$<\Psi(\th),
\Phi(\th)>=\bar\Psi^T(0)\Phi(0)-\int_{-\tau_0}^0\int_{\xi=0}^\theta\bar\Psi^T(\xi-\theta)d\eta(\theta)\phi(\xi)d\xi,$$
where $\eta(\theta)=\eta(\theta,0)$.

We assume that $\pm i\om_0$ are eigenvalues of $\ca(0)$. Thus,
they are also eigenvalues of $\ca^*$. We can easily obtain:
$$\Phi(\th)=ve^{\lam_1\th},\quad \th\in[-\tau_0, 0]\eqno(21)$$
where $v=(v_1, v_2)^T$,
$$v_1=-a\rho_{01}, v_2=\lam_1+b+a\rho_{10},$$
is the eigenvector of $\ca(0)$ corresponding to $\lam_1=i\om_0$
and
$$\Psi(s)=we^{\lam_2s},\quad s\in[0,\infty)$$ where
$w=(w_1, w_2)$,
$$w_1=\ds\f{w_2c\gamma_{10}e^{-\lam_1\tau}}{b+\lam_1+a\rho_{10}}, w_2=\ds\f{1}{\bar\eta},$$
$$\eta=v_1\ds\f{c\gamma_{10}e^{-\lam_2\tau}}{b+\lam_2+a\rho_{10}}+v_2-
\ds\f{c\gamma_{10}v_1+c\gamma_{01}v_2}{\lam_1^2}(-\tau_0\lam_1e^{-\lam_1\tau_0}-1+e^{-\lam_1\tau_0})]$$
 is the eigenvector of $\ca(0)$ corresponding to
$\lam_2=-i\om_0.$

We can verify that: $<\Psi(s), \Phi(s)>=1$, $<\Psi(s),
\bar\Phi(s)>=<\bar\Psi(s), \Phi(s)>=0$, $<\bar\Psi(s),
\bar\Phi(s)>=1.$

Using the approach of Hassard [5], we next compute the coordinates
to describe the center manifold $\Omega_0$ at $\mu=0$. Let
$u_t=u_t(t+\th), \th\in[-\tau_0,0)$ be the solution of equation
(20) when $\mu=0$ and
$$z(t)=<\Psi, u_t>,
\quad w(t,\th)=u_t(\th)-2Re\{z(t)\Phi(\th)\}.$$

On the center manifold $\Omega_0$, we have:
$$w(t,\th)=w(z(t), \bar z(t), \th)$$ where
$$w(z,\bar z, \th)=w_{20}(\th)\ds\f{z^2}{2}+w_{11}(\th)z\bar
z+w_{02}(\th)\ds\f{\bar z^2}{2}+w_{30}(\th)\ds\f{z^3}{6}+\dots$$
in which $z$ and $\bar z$ are local coordinates for the center
manifold $\Omega_0$ in the direction of $\Psi$ and $\bar\Psi$ and
$w_{02}(\th)=\bar w_{20}(\th)$. Note that $w$ and $u_t$ are real.

For solution $u_t\in \Omega_0$ of equation (20), as long as
$\mu=0$, we have:
$$\dot z(t)=\lam_1z(t)+g(z, \bar z)\eqno(22)$$ where
$$\begin{array}{ll}
g(z, \bar z)& =\bar\Psi(0)F(w(z(t),\bar z(t), 0)+2Re(z(t)\Phi(0)))=\\
& =g_{20}\ds\f{z(t)^2}{2}+g_{11}z(t)\bar z(t)+g_{02}\ds\f{\bar
z(t)^2}{2}+g_{21}\ds\f{z(t)^2\bar z(t)}{2}+\dots\end{array}$$
where
$$g_{20}=F^1_{20}\bar w_1+F^2_{20}\bar w_2,
g_{11}=F^1_{11}\bar w_1+F^2_{11}\bar w_2, g_{02}=F^1_{02}\bar
w_1+F^2_{02}\bar w_2,\eqno(23)$$ with
$$\begin{array}{l}F_{20}^1=-a(\rho_{20}v_1^2+2\rho_{11}v_1v_2+\rho_{02}v_2^2,\\
F_{20}^2=c[\gamma_{20}v_1^2e^{-2\lam_1\tau}+\gamma_{02}v_2^2e^{-2\lam_1\tau}+2\gamma_{11}v_1v_2e^{-2\lam_1\tau}],\\
F_{11}^1=-a[\rho_{20}v_1\bar v_1+\rho_{11}(v_1\bar v_2+\bar v_1v_2)+\rho_{02}v_2\bar v_2],\\
F_{11}^2=c[\gamma_{20}v_1\bar v_1+\gamma_{11}(v_1\bar v_2+\bar
v_1v_2)+\gamma_{02}v_2\bar v_2],\\
F_{02}^1= \bar F_{20}^1, F_{02}^2=\bar F_{20}^2,
\end{array}$$ and
$$g_{21}=F^1_{21}\bar w_1+F^2_{21}\bar w_2\eqno(24)$$ where
$$\begin{array}{l}
F_{21}^1= -a[\rho_{20}(2v_1w_{11}^1(0)+\bar
v_1w_{20}^1(0))+2\rho_{11}(v_1w_{11}^2(0)+\\
\ds\f{\bar v_1w_{20}^2(0)}{2}+\ds\f{\bar
v_2w_{20}^1(0)}{2}+v_2w_{11}^1(0))+ \rho_{02}(2v_2w_{11}^2(0)+\bar
v_2w_{20}^2(0))-\\\rho_{30}v_1^2\bar v_1+2\rho_{21}v_1v_2\bar
v_1+2\rho_{12}v_1v_2\bar v_2+ \rho_{03}v_2^2\bar
v_2+\rho_{21}v_1^2\bar v_2+\rho_{12}\bar v_1v_2^2]\end{array}$$
$$\begin{array}{l}F_{21}^2=c[\gamma_{20}(2v_1w_{11}^1(-\tau)e^{-\lam_1\tau}+\bar v_1w_{20}^1(-\tau)e^{\lam_1\tau})
+2\gamma_{11}(v_1w_{11}^2(-\tau)e^{-\lam_1\tau}+\\+ \ds\f{\bar
v_1w_{20}^2(-\tau)e^{\lam_1\tau}}{2}\!+\!\ds\f{\bar
v_2w_{20}^1(-\tau)e^{\lam_1\tau}}{2}\!+\!v_2w_{11}^1(-\tau)e^{-\lam_1\tau})
\!+\!\gamma_{02}(2v_2w_{11}^2(-\tau)e^{-\lam_1\tau}+\\+\bar
v_2w_{20}^2(-\tau)e^{\lam_1\tau}) +\gamma_{30}v_1^2\bar
v_1e^{-\lam_1\tau}+\gamma_{21}(2v_1\bar
v_1v_2e^{-\lam_1\tau}+v_1^2\bar v_2e^{-\lam_1\tau})+\\
+\gamma_{12}(2v_1v_2\bar v_2e^{-\lam_1\tau}+\bar
v_1v_2^2e^{-\lam_1\tau})+\gamma_{03}v_2^2\bar
v_2e^{-\lam_1\tau}].\end{array}$$

The vectors $w_{20}(\theta)$, $w_{11}(\theta)$ with
$\theta\in[-\tau,0]$ are given by:
$$\begin{array}{l}
w_{20}(\theta)=-\ds\f{g_{20}}{\lam_1}ve^{\lam_1\theta}-\ds\f{\bar
g_{02}}{3\lam_1}\bar ve^{\lam_2\theta}+E_1e^{2\lam_1\theta}\\
w_{11}(\theta)=\ds\f{g_{11}}{\lam_1}ve^{\lam_1\theta}-\ds\f{\bar
g_{11}}{\lam_1}\bar ve^{\lam_2\theta}+E_2\end{array}\eqno(25)$$
where
$$
E_1=-(A+e^{-2\lam_1\tau_0}B-2\lam_1I)^{-1}F_{20},\quad
E_2=-(A+B)^{-1}F_{11},$$ where $F_{20}=(F_{20}^1, F_{20}^2)^T$,
$F_{11}=(F_{11}^1, F_{11}^2)^T$.

Based on the above analysis and calculation, we can see that each
$g_{ij}$ in (23), (24) are determined by the parameters and delay
of system (2). Thus, we can explicitly compute the following
quantities:
$$\begin{array}{l}
C_1(0)=\ds\f{i}{2\om_0}(g_{20}g_{11}-2|g_{11}|^2-\ds\f{1}{3}|g_{02}|^2)+\ds\f{g_{21}}{2},\\
\mu_2=-\ds\f{Re(C_1(0))}{M}, T_2=-\ds\f{Im(C_1(0))+\mu_2N}{\om_0},
\beta_2=2Re(C_1(0)),\end{array}\eqno(26)$$ where $M$ and $N$ are
given by (16) and (17).
%$$\lam'=-\ds\f{e^{-\lam \tau}\lam (\lam q_1+q_0)}
%{2\lam^1+p_1+e^{-\lam\tau}(-q_1+\tau(q_1\lam q_1+q_0))}$$

In short, this leads to the following result:

{\bf Theorem 3.} {\it In formulas (26), $\mu_2$ determines the
direction of the Hopf bifurcation: if $\mu_2>0 (<0)$, then the
Hopf bifurcation is supercritical (subcritical) and the
bifurcating periodic solutions exist for $\tau>\tau_0 (<\tau_0)$;
$\beta_2$ determines the stability of the bifurcating periodic
solutions: the solutions are orbitally stable (unstable) if
$\beta_2<0 (>0)$; and $T_2$ determines the period of the
bifurcating periodic solutions: the period increases (decreases)
if $T_2>0 (<0)$.}

\vspace{0.6cm}
\section*{\normalsize\bf 4. Numerical example.}

\hspace{0.6cm} In this section we find the waveform plots through
the formula:
$$X(t+\th)\!=\! z(t)\Phi(\th)\!+\!\bar
z(t)\bar\Phi(\th)\!+\!\ds\f{1}{2}w_{20}(\th)z^2(t)+w_{11}(\th)z(t)\bar
z(t)\!+\!\ds\f{1}{2}w_{02}(\th)\bar z(t)^2+X_0,$$ where $z(t)$ is
the solution of (20), $\Phi(\th)$ is given by (21), $w_{20}(\th),
w_{11}(\th), w_{02}(\th)=\bar w_{20}(\th)$ are given by (25) and
$X_0=(y_{10}, y_{20})^T$ is the equilibrium state.

For the numerical simulations we use Maple 9.5. and the data from
[10]: the degradation of p53 through ubiquintin pathway
$a=3\times10^{-2}$sec${^-1}$, the spontaneous degradation of P53
is $b=10^{-4}$sec$^{-1}$, the dissociation constant between P53
and Mdm2 gene is $k_1=28$, the degradation rate of Mdm2 protein is
$d=10^{-2}$sec$^{-1}$, the P53 protein production rate is
$S=0.01$sec$^{-1}$ and the production rate of Mdm2 is
$c=1$sec$^{-1}$. For this date we consider the different values
for the constant k. By changing instantaneously the dissociation
constant k the response of the system is different.

For $k=17.5$ we obtain: the equilibrium point
$y_{10}=1.674122637$, $y_{20}=4.587857801$, the coefficients which
describe the limit cycle: $\mu_2\!=\! 0.002340098338$,
$\beta_2\!=\!-0.000001452701020$, $T_2\!=\!0.00007793850775$ and
$\omega_0\!=\! 0.04577286901$, $\tau_0\!=\! 60.52296388$. Then the
Hopf bifurcation is supercri\-ti\-cal, the solutions are orbitally
stable and the period of the solution is increasing. The wave
plots are displayed in fig1 and fig2:

\begin{center}
{\small \begin{tabular}{c|c} \hline Fig.1. $(t,y_1(t))$&Fig.2.
$(t,y_2(t))$\\&\\
 \cline{1-2} \epsfxsize=5cm

\epsfysize=5cm

\epsffile{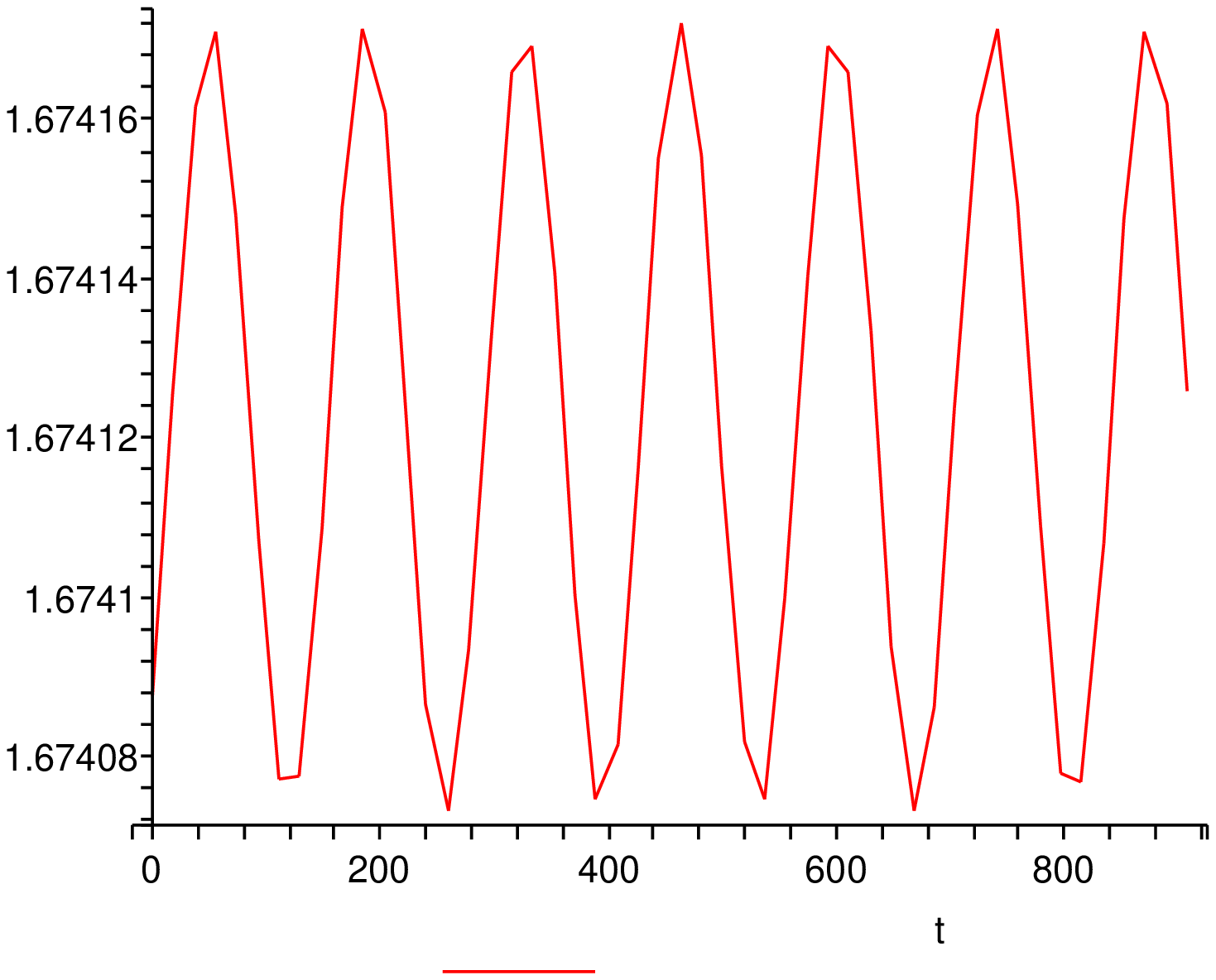} &

\epsfxsize=6cm

\epsfysize=5cm

\epsffile{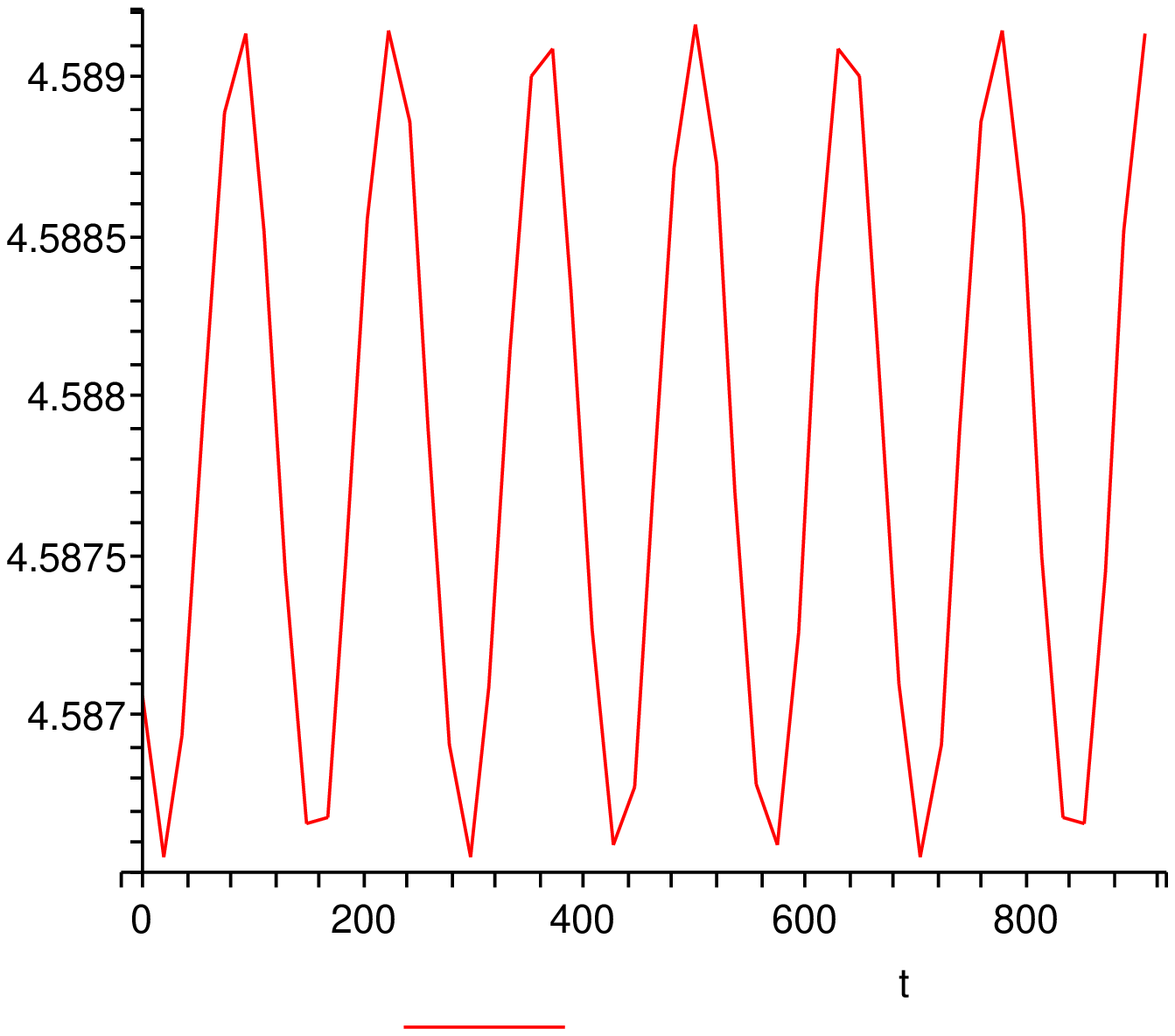}

\\
 \hline
\end{tabular}}
\end{center}

\medskip

For $k=120$ we obtain: the equilibrium point $y_{10}=3.853801769$,
$y_{20}=11.20502079$, the coefficients which describe the limit
cycle: $\mu_2\!=\!-8.219839362\cdot 10^{-7}$,
$\beta_2\!=\!2.0434666644\cdot 10^{-10}$,
$T_2\!=\!2.243699328\cdot 10^{-7}$ and $\omega_0\!=\!
0.02923152517$, $\tau_0\!=\! 92.68683554$. Then the Hopf
bifurcation is subcri\-ti\-cal, the solutions are orbitally
unstable and the period of the solution is increasing. The wave
plots are presented in fig3 and fig4:

\begin{center}
{\small \begin{tabular}{c|c} \hline Fig.3. $(t,y_1(t))$&Fig.4.
$(t,y_2(t))$\\&\\
 \cline{1-2} \epsfxsize=5cm

\epsfysize=5cm

\epsffile{fig1.eps} &

\epsfxsize=6cm

\epsfysize=5cm

\epsffile{fig2.eps}

\\
 \hline
\end{tabular}}
\end{center}

\medskip

For $k=1750$ we obtain: the equilibrium point
$y_{10}=14.88816840$, $y_{20}=34.27918463$, the coefficients which
describe the limit cycle: $\mu_2\!=\!1.104273140\cdot 10^{-11}$,
$\beta_2\!=\!-6.4139104\cdot 10^{-16}$, $T_2\!=\!1.924908056\cdot
10^{-11}$ and $\omega_0\!=\! 0.01423761906$, $\tau_0\!=\!
174.4149631$. Then the Hopf bifurcation is supercri\-ti\-cal, the
solutions are orbitally stable and the period of the solution is
increasing. The wave plots are given in fig5 and fig6:

\begin{center}
{\small \begin{tabular}{c|c} \hline Fig.5. $(t,y_1(t))$&Fig.6.
$(t,y_2(t))$\\&\\
 \cline{1-2} \epsfxsize=5cm

\epsfysize=5cm

\epsffile{fig1.eps} &

\epsfxsize=6cm

\epsfysize=5cm

\epsffile{fig2.eps}

\\
 \hline
\end{tabular}}
\end{center}

\medskip

\section*{\normalsize\bf 5. Conclusions.}

\hspace{0.6cm} For the present model, we obtain an oscillatory
behavior similar with the findings in [10], according with the
qualitative study.

We have proved that a limit cycle exists and it is characterized
by the coefficients from (26). For different values for the
equilibrium constant k, in Section 4 we obtain stable or unstable
periodic solutions with increasing periods, via a Hopf
bifurcation.


\begin{thebibliography}{99}

\bibitem{}  Adimy M., Crauste F., Halanay  A., Neam\c tu M., Opri\c s D., Stability of limit cycle in a
pluripotent stem cell dynamics model, Chaos, Solitons and Fractals
J., 27(2006), 1091-1107.
\bibitem{} Chickarmane V., Nadim A., Ray A., Sauro H.M., A P53
oscillator model of DNA break repair control,
arXiv:q-bio.MN/0510002v1.
\bibitem{}   Lahav  G., Rosenfeld N., Sigal A., Geva-Zatorsky  N., Levine A.J., Elowitz M.B., Alon U.,
Dynamics of the p53-Mdm2 feedback loop in individual cells, Nat.
Genet., 36(2004), 147-150.
\bibitem{}  Lev Bar-Or R., Maya R., Segel L.A., Alon U., Levine A.J., Oren M., Generation of oscillations by p53-Mdm2 feedback loop:
A theoretical and experimental study, PNAS, 97(2000), no.21,
11250-11255.
\bibitem{} Hassard B.D., Kazarinoff N.D., Wan Y.H., Theory and applications of Hopf bifurcation, Cambridge University Press, Cambridge,
1981.
\bibitem{} Kohn K. W., Pommier Y., Molecular interaction map of p53 and Mdm2 logic elements,
which control the Off-On switch of p53 in response to DNA damage,
Science Direct, Biochemical and Biophysical Research
Communications, 331(2005), 816-827.
\bibitem{}  Mihala\c s G.I., Neam\c tu M., Opri\c s D.,
Horhat R.F., A dynamic P53-Mdm2 model with time delay, it will be
appear in Chaos, Solitons and Fractals J.
\bibitem{}  Mihala\c s G.I., Simon Z., Balea  G., Popa E., Possible oscillatory
behaviour in p53-Mdm2 interaction computer simulation, J. of
Biological Systems, 8(2000), nr. 1, 21-29.
\bibitem{} Perry M.E., Mdm2 in the response to radiation, Mol. Cancer Res., 2(2004),
9-19.
\bibitem{} Tiana G., Jensen M.H., Sneppen K., Time delay as a key
to apoptosis induction in the P53 network, Eur. Phys. J.,
29(2002), 135-140.

\end{thebibliography}
\end{document}